\documentclass[12pt]{amsart}

\usepackage{amsmath} \usepackage{amssymb} \usepackage{amsthm}
\usepackage{amscd}

\theoremstyle{plain}
\newtheorem{thm}{Theorem}[section]

\newtheorem{prop}[thm]{Proposition}

\theoremstyle{definition}
\newtheorem{defn}[thm]{Definition}
\newtheorem{exam}[thm]{Example}
\theoremstyle{remark}
\newtheorem{remark}[thm]{Remark}
\newtheorem{quest}[thm]{Question}
\DeclareMathOperator{\codim}{codim}

\DeclareMathOperator{\fpt}{fpt} \DeclareMathOperator{\lct}{lct}
\newcommand{\frm}{{\mathfrak m}}

\newcommand{\fra}{{\mathfrak a}}
\newcommand{\frb}{{\mathfrak b}}

\newcommand{\frq}{{\mathfrak q}}

\newcommand{\bbR}{\ensuremath{\mathbb R}}
\newcommand{\bbZ}{\ensuremath{\mathbb Z}}

\newcommand{\bbQ}{\ensuremath{\mathbb Q}}
\newcommand{\bbF}{\ensuremath{\mathbb F}}

\newcommand{\JJ}{\ensuremath{\mathcal J}}

\newcommand{\OO}{\ensuremath{\mathcal O}}

\newcommand{\llbracket}{[\negthinspace[}
\newcommand{\rrbracket}{]\negthinspace]}

\title{Test ideals vs. multiplier ideals}
\author{Mircea Musta\c t\u a}
\address{Department of Mathematics, University of Michigan, Ann Arbor, MI 48109, USA}
\email{mmustata@umich.edu}
\author{Ken-ichi Yoshida}
\address{Graduate School of Mathematics, Nagoya University, Chikusa-ku,
Nagoya,464-8602 JAPAN} \email{yoshida@math.nagoya-u.ac.jp}
\subjclass[2000]{Primary: 13A35; Secondary: 14B05}

\dedicatory{}
\date{\today}

\begin{document}

\begin{abstract}
The generalized test ideals introduced in \cite{HY} are related to
multiplier ideals via reduction to characteristic $p$. In addition,
they satisfy many of the subtle properties of the multiplier ideals,
which in characteristic zero follow via vanishing theorems. In this
note we give several examples to emphasize the different behavior of
test ideals and multiplier ideals. Our main result is that every
ideal in an $F$-finite regular local ring can be written as a
generalized test ideal. We also prove the semicontinuity of $F$-pure
thresholds (though the analogue of the Generic Restriction Theorem
for multiplier ideals does not hold).
\end{abstract}

\thanks{The first author was partially supported
 by NSF grant DMS 0500127 and by a Packard Fellowship.}

\maketitle

\markboth{M.~MUSTA\c T\u A and K.~YOSHIDA}{TEST IDEALS VS.
MULTIPLIER IDEALS}

\section{Introduction}

In recent years the multiplier ideals and the log canonical
threshold have played an important role in higher dimensional
birational geometry (see e.g. \cite{Laz}). These are invariants of
singularities in characteristic zero, that can be defined in terms
of log resolutions of singularities. Suppose for simplicity that $X$
is a smooth variety over a field of characteristic zero, and that
$\fra\subseteq\OO_X$ is a coherent sheaf of ideals. The multiplier
ideal associated to the pair $(X,\fra)$ and to a non-negative real
number $t$ is denoted by $\JJ(\fra^t)$. If $t_1<t_2$, then
$\JJ(\fra^{t_2})\subseteq\JJ(\fra^{t_1})$, and $\JJ(\fra^t)=\OO_X$
for $0<t\ll 1$. The smallest $t$ such that $\JJ(\fra^t)\neq\OO_X$ is
the \emph{log canonical threshold} $\lct(\fra)$.

On the other hand, in positive characteristic one can define
invariants using the Frobenius morphism. Specifically, Hara and the
second author introduced in \cite{HY} a notion of tight closure for
pairs, and corresponding (generalized) test ideals $\tau(\fra^t)$.
Suppose that we have a pair $(X,\fra)$ and $t\in\bbR_+$, where $X$
is a smooth variety over a field of characteristic zero. If we
denote by $\frb_p$ the reduction mod $p$ of the ideal $\frb$, it was
proved in \cite{HY} that
\begin{equation}\label{eq0}
\JJ(\fra^t)_p=\tau(\fra_p^t)
\end{equation}
 for all primes $p\gg 0$ (depending on
$t$).

\par
In the same vein,
 Takagi and Watanabe defined in
 positive characteristic \cite{TW} the \emph{F-pure threshold}
$\fpt(\fra)$. When the ambient variety is nonsingular and $F$-finite
(that is, the Frobenius morphism $F\colon X\to X$ is finite), this
can be described as the smallest $t$ such that
$\tau(\fra^t)\neq\OO_X$. The formula (\ref{eq0}) can then be
reinterpreted as saying that
\begin{equation}\label{eq00}
\lim_{p\to\infty}\fpt(\fra_p)=\lct(\fra).
\end{equation}

The above shows the close connection between multiplier and test
ideals. In fact, more is true. Multiplier ideals satisfy several
subtle properties, such as the Restriction Theorem, the
Subadditivity and the Summation Theorems, and Skoda's Theorem (see
\cite{Laz}). One common feature of these results is that they all
rely on applications of vanishing theorems. As it was pointed out in
\cite{HY}, \cite{HT} and \cite{Ta}, all these results have similar
statements for test ideals, with substantially easier proofs.

On the other hand, multiplier ideals enjoy several other properties,
that follow simply from the description in terms of resolutions of
singularities. In this note we concentrate on these properties, and
show that essentially all these fail for test ideals.

\bigskip

Our basic ingredient is the description of test ideals from
\cite{BMS1}, which holds when the ambient variety is nonsingular and
$F$-finite. Therefore we will always make this assumption. Our main
result is a positive one: under mild assumptions, every ideal is a
test ideal.

\begin{thm}\label{thm1}
Suppose that $R$ is a ring of characteristic $p>0$, such that $R$ is
a finitely generated free module over $R^p$. For every ideal $I$ in
$R$, there is $f\in R$ and $c>0$ such that $I=\tau(f^c)$.
\end{thm}

Note that the theorem applies when $R$ is a local regular $F$-finite
ring, or when $R=k[x_1,\ldots,x_n]$, where $[k\colon k^p]<\infty$.
 As we will see, both $f$ and $c$ in the theorem can be
explicitly determined. Moreover, if $I$ is $\frm$-primary, for some
maximal ideal $\frm$, then we show that we may write also
$I=\tau(\fra^{c'})$ for some $\frm$-primary ideal $\fra$ and some
$c'>0$.

Note that Theorem~\ref{thm1} contrasts with the situation for
multiplier ideals. In that case, as an immediate consequence of the
definition one shows that every multiplier ideal is integrally
closed. Moreover, as it was recently shown in \cite{LL}, there are
more subtle conditions involving the local syzygies, that are
satisfied by all multiplier ideals.

In \cite{ELSV} one shows that whenever one writes an ideal $I$ as a
multiplier ideal, then one can prove an effective uniform Artin-Rees
theorem for $I$. The main ingredient in that proof is a basic
property of multiplier ideals that follows from the definition via
resolutions. As we show in Example~\ref{AR} below, this property
fails in the case of test ideals, and therefore it seems that
Theorem~\ref{thm1} does not have similar consequences in the
direction of uniform Artin-Rees statements.

We give several examples to illustrate that basic properties of
multiplier ideals, which easily follow from the definition via log
resolutions, can fail in the case of test ideals:
\begin{enumerate}
\item[i)] We show that it can happen that for a (principal) ideal
$\fra$, we can have the ideal $\tau(\fra^c)$ non-radical, where
$c=\fpt(\fra)$ (see Example~\ref{counterex}).
\item[ii)] We give an example of a (principal) ideal $\fra$ with $c=\fpt(\fra)$
such that $\tau(\fra^c)$ is $\frm$-primary for a maximal ideal
$\frm$, but such that $\fpt(\fra)<\fpt(\fra+\frm^{\ell})$ for all
$\ell\gg 0$ (see Example~\ref{counterex}).
\item[iii)] We show that the analogue of
the Generic Restriction Theorem for multiplier ideals can fail in
the case of test ideals (see Example~\ref{semicont1}). However, we
will prove that the $F$-pure thresholds satisfy the same
semicontinuity property as the log canonical thresholds.
\end{enumerate}

The paper is structured as follows. In \S 2 we review the
definitions of multiplier and generalized test ideals, and some
basic properties. In particular, we recall the description of test
ideals in the case of a regular $F$-finite ring from \cite{BMS1},
which we will systematically use. In \S 3 we prove
Theorem~\ref{thm1} above. The next section is devoted to various
examples, including the ones mentioned above, while in the last
section we prove the semicontinuity result for $F$-pure thresholds.

\section{Preliminaries}

We first recall the definition of multiplier ideals (for details see
\cite{Laz}, \S 9). For a real number $u$, we denote by $\lceil
u\rceil $ the smallest integer $\geq u$. Similarly, $\lfloor
u\rfloor$ is the largest integer $\leq u$. This notation is extended
to divisors with real coefficients, in which case we apply it to
each coefficient.

Let $X$ be a $\bbQ$-Gorenstein normal variety over a field of
characteristic zero, $Y \subsetneq X$ a proper closed subscheme
defined by an ideal sheaf $\fra \subseteq \OO_X$, and $t \ge 0$ a
real number. Suppose that $\pi \colon \widetilde{X} \to X$ is a log
resolution of the pair $(X,Y)$ such that $\fra\OO_{\widetilde{X}} =
\OO_{\widetilde{X}}(-F)$, and let $K_{\widetilde{X}/X}$ denote the
discrepancy divisor. Then the \textit{multiplier ideal}
$\JJ(\fra^t)$ is defined by
\[
 \JJ(\fra^t) = \pi_{*}\OO_{\widetilde{X}}
\bigg(\lceil K_{\widetilde{X}/X} - t F \rceil \bigg) \subseteq
\OO_{X}.
\]
This is an ideal of $\OO_X$ that does not depend on the choice of
the log resolution.

One says that $X$ has \emph{log terminal singularities} at $x\in X$
if $x$ does not lie in the support of $\JJ(\fra^t)$ for $0<t\ll 1$.
In this case one defines the \textit{log canonical threshold} of
$\fra$ at $x$, denoted by $\lct_x(\fra)$, to be
\[
 \lct_x(\fra) = \sup \{ s \in \bbR_{\ge 0} \mid
x\,\text{is not in the support of}\,\JJ(\fra^s)\}.
\]

For the purpose of this paper, it is enough to restrict ourselves to
the case when the variety $X$ is nonsingular (hence, in particular,
$X$ has log terminal singularities at every point). It is easy to
see starting from definition that $\JJ(\fra^{t_1})\subseteq
\JJ(\fra^{t_2})$ if $t_1>t_2$. Moreover, given any $t\geq 0$, there
is a positive $\varepsilon$ such that
$\JJ(\fra^t)=\JJ(\fra^{t+\varepsilon})$. Following \cite{ELSV}, we
call $\lambda>0$ a \emph{jumping number} of $\fra$ if
$\JJ(\fra^{\lambda})\neq\JJ(\fra^{t})$ for every $t<\lambda$. With
the notation in the definition of multiplier ideals, it follows
easily that if we write $F=\sum_ia_iE_i$, then for every jumping
number $\lambda$ of $\fra$, there is $i$ such that $a_i\lambda$ is
an integer. In particular, the jumping numbers are rational and they
form a discrete set.

The smallest jumping number of $\fra$ is the \emph{log canonical
threshold} $\lct(\fra)$. It is clear that we can define local
versions of the jumping numbers at every $x\in X$. In this case, the
smallest jumping number is precisely $\lct_x(\fra)$. In fact, it is
easy to see that $\lct(\fra)=\min_{x\in X}\lct_x(\fra)$.

\bigskip

We now turn to the positive characteristic setting. Let $R$ be a
Noetherian ring containing a field of characteristic $p>0$. The ring
$R$ is called \textit{$F$-finite} if $R$ is a finitely generated
module over its subring $R^p = \{a^p \in R\,:\, a \in R\}$. If $J$
is an ideal in $R$, then $J^{[p^e]}$ denotes the ideal
$(u^{p^e}\colon u\in J)$.
 We recall first the notion of \textit{generalized test ideals},
introduced by Hara and the second author in \cite{HY}. We denote by
$R^{\circ}$ the complement of all minimal prime ideals of $R$.

\begin{defn} \label{hy-tau}
Let $\fra$ be an ideal such that $\fra \cap R^{\circ} \ne
\emptyset$. Let $t \ge 0$ be a real number. For any ideal $I$ of
$R$, the \textit{$\fra^t$-tight closure} of $I$, denoted by
$I^{*\fra^t}$, is defined to be the ideal of $R$ consisting of all
elements $z \in R$ for which
 there exists $c \in R^{\circ}$ such that
\[
 cz^q \fra^{\lceil t q \rceil}  \subseteq I^{[q]}
\]
for all large $q=p^e$. Here we denote by $\lceil x\rceil$ the
smallest integer $\geq x$.

\par
Assume that $R$ is excellent and reduced. Given a real number $t \ge
0$, one defines the \emph{generalized test ideal} $\tau(\fra^t)$ by
\[
 \tau(\fra^t) = \bigcap_{I \subseteq R} I \colon I^{*\fra^t},
\]
where $I$ runs through all ideals of $R$. In the case of a principal
ideal $\fra=(f)$, we simply write $\tau(f^t)$.
\end{defn}

\par
Blickle, Smith and the first author gave in \cite{BMS1} a different
description of generalized test ideals in the case of an $F$-finite
regular ring $R$. We briefly recall this description here, in the
special case when $R$ is free and finitely generated over $R^p$.
Note that this condition holds, for example, when $R$ is an
$F$-finite regular local ring, or when $R=k[x_1,\ldots,x_n]$ and
$[k\colon k^p]<\infty$.

It follows from our assumption that for every $p^e$, with $e\geq 1$,
$R$ is free over $R^{p^e}=\{a^{p^e}\colon a\in R\}$. For every such
$e$, let us fix a basis $u_1,\ldots,u_N$ of $R$ over $R^{p^e}$.
Given any ideal $\frb$ of $R$, we choose generators $h_1,\ldots,h_s$
of $\frb$. If we write for every $i$
\[
 h_i = \sum_{j=1}^N a_{ij}^{p^e}u_j,
\]
with $a_{ij} \in R$, then we put
\[
 \frb^{[1/p^e]} = (a_{ij}\colon 1\le i \le s,\, 1 \le j \le N).
\]
In fact, $\frb^{[1/p^e]}$ is the unique smallest ideal $J$ (with
respect to inclusion) such that $\frb\subseteq J^{[p^e]}$. In
particular, $\frb^{[1/p^e]}$ does not depend on the choice of
generators for $\frb$, or on the choice of basis for $R$ over
$R^{p^e}$.

Suppose now that $\fra$ is an ideal in $R$ and that $t$ is a
positive real number. For every $e\geq 1$ we have the inclusion
\[
\left(\fra^{\lceil tp^e\rceil}\right)^{[1/p^e]} \subseteq
\left(\fra^{\lceil tp^{e+1}\rceil}\right)^{[1/p^{e+1}]}.
\]
Since $R$ is Noetherian, these ideals stabilize for $e\gg 0$, and
the limit was taken as definition for $\tau(\fra^t)$ in \emph{loc.
cit}, the equivalence with the definition from \cite{HY} being
proved in \emph{ibid.}, Proposition~2.22.

\bigskip

We now recall the definition of $F$-jumping exponents, that is
analogous to that of jumping  numbers for multiplier ideals. We
assume that $R$ is a regular $F$-finite ring. Note that if $t < t'$,
then $\tau(\fra^t) \supseteq \tau(\fra^{t'})$. Moreover, for every
$t$ there exists $\varepsilon >0$ such that
$\tau(\fra^{t})=\tau(\fra^{t'})$ for every $t' \in
[t,t+\varepsilon)$.

\begin{defn} \label{f-jump}
A positive real number $\lambda$ is called an \textit{$F$-jumping
exponent} of $\fra$ if $\tau(\fra^\lambda) \ne \tau(\fra^{t})$ for
every $t < \lambda$. It is convenient to make also the convention
that $0$ is an $F$-jumping exponent.
\end{defn}

The smallest positive $F$-jumping exponent of $\fra$ is the
\emph{$F$-pure threshold} $\fpt(\fra)$. This notion was introduced
in a more general setting by Takagi and Watanabe in \cite{TW}, as an
analogue of the log canonical threshold.

When $(R,\frm)$ is an $F$-finite regular local ring, the $F$-pure
threshold has the following alternative description (see \cite{BMS1}
or \cite{MTW}). Given an ideal $\fra\subseteq\frm$ and $e\geq 1$, we
denote by $\nu(e)$ the largest integer $r$ such that
$\fra^r\not\subseteq\frm^{[p^e]}$ (we put $\nu(e)=0$ if there is no
such $r$). We then have
\begin{equation}\label{fpt}
\fpt(\fra)=\sup_e\frac{\nu(e)}{p^e}.
\end{equation}
It follows that given a nonnegative integer $c$, we have
$\fpt(\fra)\leq c$ if and only if $\fra^{\lfloor
cp^e\rfloor+1}\subseteq \frm^{[p^e]}$ for every $e$.

\bigskip

Rationality and discreteness of $F$-jumping exponents is more subtle
in positive characteristic. Both properties have been proved in
\cite{BMS1} for an arbitrary ideal in a regular ring that is
essentially of finite type over an $F$-finite field, and for a
principal ideal in any $F$-finite regular ring in \cite{BMS2}.

We will be especially interested in the case when $\fra=(f)$ is a
principal ideal in an $F$-finite regular ring. In this case, Skoda's
Theorem (see Theorem~4.1 in \cite{HT} or Proposition~2.25 in
\cite{BMS1}) implies that for every $t \ge 1$ we have $\tau(f^t) = f
\cdot \tau(f^{t-1})$. Therefore the set of $F$-jumping exponents of
$f$ is periodic with period one, hence it is enough to describe the
$F$-jumping exponents in the interval $(0,1]$. As we have mentioned,
this is a finite set.

\section{Any ideal in an $F$-finite regular local ring is a test ideal}

Throughout this section we assume that $R$ is a regular, $F$-finite
ring. By a theorem of Kunz \cite{Kunz}, this is equivalent with $R$
being finitely generated and projective over $R^p$. We will assume
that moreover, $R$ is free over $R^p$. This holds, for example, if
$R$ is also local, or if $R=k[x_1,\ldots,x_n]$, where $[k\colon
k^p]<\infty$. The following is the main result of this section.

\begin{thm}\label{thm2}
Let $R$ be a regular ring of characteristic $p>0$, such that $R$ is
a finitely generated, free module over $R^p$.
\begin{enumerate}
\item[1)] For every ideal $I$ in $R$, there are $f\in R$ and $c>0$ such
that $I=\tau(f^c)$.
\item[2)] Moreover, if $\frm$ is a maximal ideal in $R$, and if
$I$ is $\frm$-primary, then we can find an $\frm$-primary ideal
$\frb$ and $c'>0$ such that $I=\tau(\frb^{c'})$.
\end{enumerate}
\end{thm}

Suppose that $R$ satisfies the hypothesis of the theorem, and let
$N={\rm rk}_{R^p}(R)$. It is clear that $N=1$ if and only if
$\dim(R)=0$, in which case Theorem~\ref{thm2} is trivial. We will
henceforth assume $N>1$. Note that if $e\geq 1$, then $R$ is free
over $R^e$ of rank $N^e$.

The first assertion in Theorem~\ref{thm2} follows from the more
precise statement below.

\begin{prop} \label{Precise}
Let $R$ be a ring of characteristic $p>0$ that is free and finitely
generated over $R^p$, with ${\rm rk}_{R^p}(R)=N$. Let
$I=(z_1,\ldots,z_{\mu})$ be an ideal of $R$, and fix $e_0\geq 1$
such that $N^{e_0}\geq\mu$. If $g_1,\ldots,g_{N^{e_0}}$ is a basis
of $R$ over $R^{p^{e_0}}$, and if we put
\[
 f=\sum_{i=1}^{\mu} z_i^{p^{e_0}} g_i \in R, \qquad c=\frac{1}{p^{e_0}} \in \bbQ,
\]
then
\[
 \tau(f^c)=I.
\]
\end{prop}

\begin{proof}
We use the description of $\tau(f^c)$ from \cite{BMS1}. If $e\geq
e_0$, then we have a basis of $R$ over $R^{p^e}$ given by
\[
\{g_{i_1}g_{i_2}^p \cdots g_{i_{e-e_0+1}}^{p^{e-e_0}}\mid 1\leq
i_1,\ldots,i_{e-e_0+1}\leq N\}.
\]
Since we can write
$f^{p^{e-e_0}}=\sum_{i=1}^{\mu}z_i^{p^e}g_i^{p^{e-e_0}}$, it follows
that
\[
\left(f^{\lceil cp^e\rceil}\right)^{[1/p^e]}=
\left(f^{p^{e-e_0}}\right)^{[1/p^e]}=(z_1,\ldots,z_{\mu})=I.
\]
Since this is true for every $e\geq e_0$, we deduce $\tau(f^c)=I$.
\end{proof}

\par \vspace{2mm}
We turn now to the second assertion in Theorem~\ref{thm2} (this
answers positively a question raised by Kei-ichi Watanabe). The
assertion is a consequence of 1), together with the more general
statement below. Recall that by Corollary~2.16 in \cite{BMS1}, for
every $f$ and every $c$ there is $\varepsilon>0$ such that
$\tau(f^c)=\tau(f^{c+\varepsilon})$.

\begin{prop} \label{mprimary}
Let $R$ be a regular $F$-finite ring, and $\frm$ a maximal ideal in
$R$. Suppose that $f\in R$ and $c>0$ are such that $I:=\tau(f^c)$ is
$\frm$-primary. If we fix $\varepsilon
> 0$ such that $I=\tau(f^{c+\varepsilon})$, and if $r$ is such that
$\frm^r\subseteq I$, then for every positive integer $\ell$ with
$\ell \varepsilon \ge r+\codim(\frm)-1$, we have
\[
 I=\tau((fR+\frm^{\ell})^{c+\varepsilon}).
\]
\end{prop}

\begin{proof}
We put $\fra_{\ell} = fR+\frm^{\ell}$. Note that we clearly have $I
=\tau(f^{c+\varepsilon}) \subseteq
\tau(\fra_{\ell}^{c+\varepsilon})$.
\par \vspace{2mm}

On the other hand, by Takagi's Summation Theorem (see Theorem~3.1 in
\cite{Ta}), we have
\[
 \tau(\fra_{\ell}^{c+\varepsilon})
\subseteq \sum_{\lambda+\nu=c+\varepsilon}
\tau(f^{\lambda})\cdot\tau(\frm^{\ell\nu}) \subseteq \tau(f^c) +
\tau(\frm^{\ell\varepsilon}).
\]
For the second inclusion we used the fact that if $\lambda\ge c$,
then $\tau(f^{\lambda})\subseteq\tau(f^c)$, and otherwise we have
$\nu\geq\varepsilon$, hence $\tau(\frm^{\ell\nu})\subseteq
\tau(\frm^{\ell\varepsilon})$.

\par  \vspace{2mm}
Since $\ell \varepsilon \ge r + d -1$, where $d=\codim(\frm)$, and
since $\tau(\frm^{\alpha})= \frm^{\lfloor \alpha\rfloor-d+1}$ for
every $\alpha\geq d-1$, it follows that
\[
 \tau(\frm^{\ell\varepsilon})
\subseteq \frm^r\subseteq I.
\]
Therefore $\tau(\fra_{\ell}^{c+\varepsilon})\subseteq I$, which
completes the proof of the proposition.
\end{proof}

\bigskip

Let $I$ be an ideal of a ring $R$. Recall that the \textit{integral
closure} of $I$, denoted  $\overline{I}$, is the ideal of $R$
consisting of all $z$ that satisfy an equation $f(z)=0$ for some
\[
 f(X) = X^n + a_1 X^{n-1} + \cdots + a_n \quad (a_i \in I^i).
\]
The ideal $I$ is \textit{integrally closed} if $I=\overline{I}$. It
is an immediate consequence of the definition that all multiplier
ideals are integrally closed (see \cite{Laz}, Corollary~9.6.13).

In positive characteristic, the generalized test ideal of
$\tau(\fra^t)$ is integrally closed for every $t \in \bbR_{\ge 0}$
if $\fra$ is generated by monomials in a polynomial ring (in fact,
in this case, the test ideals are given by the same formula as the
multiplier ideals in characteristic zero, see Theorem~6.10 in
\cite{HY}). More precisely, if the ideal $\fra$ is generated by
monomials in a polynomial ring, then
\[
 \tau(\fra^t)
= \left\{x^u \in R \mid u+(1,1,\ldots,1) \in {\rm Int}(t\cdot
P(\fra)) \right\},
\]
where $P(\fra)$ is the Newton polyhedron associated to $\fra$.

We mention that in dimension two, Lipman and Watanabe \cite{LW} and
Favre and Jonsson \cite{FJ} independently proved that every
integrally closed ideal is the multiplier ideal of some ideal. There
was some belief that such a result would be true in higher
dimensions. However, recent work of Lazarsfeld and Lee \cite{LL}
shows that in fact multiplier ideals have to satisfy also some
strong properties in terms of their local syzygies, allowing to give
examples in dimension $\geq 3$ of integrally closed ideals that are
not multiplier ideals.

However, as Theorem~\ref{thm2} clearly shows, the situation for test
ideals in positive characteristic is drastically different. Since
any ideal is a test ideal, in particular we get many non-integrally
closed test ideals. Here is a concrete such example.

\begin{exam}\label{taunotint}
Let $R = \bbF_2[[x,y,z]]$ and $f = x^2 +y^5+z^5$. It follows from
Proposition~\ref{Precise} that $\tau(f^{1/2})=(x,y^2,z^2)$, hence it
is \textit{not} integrally closed. In fact, we will see in
Proposition~\ref{concrete} below that $f$ has no jumping numbers in
$(1/2, 1)$. It follows that we may apply Proposition~\ref{mprimary}
with $\varepsilon=5/11$ and $r=3$ to deduce that if
$\fra=(f)+(x,y,z)^{11}$, then $\tau(\fra^{21/22})=(x,y^2,z^2)$.
\end{exam}

\begin{remark} \label{twodim}
Suppose that $(R,\frm)$ is a two-dimensional excellent Gorenstein
$F$-rational local domain of characteristic $p>0$. If $\fra
\subseteq R$ is an $\frm$-primary integrally closed ideal, and if
$\frb$ is its minimal reduction, then $\tau(\fra) = \frb : \fra$,
hence $\tau(\fra)$ is integrally closed. See \cite[Theorem 3.1]{HWY}
and \cite[Theorem 5.1]{HY}.
\end{remark}

\begin{remark} \label{poly}
In the case of a polynomial ring we do not need the assumption that
the ring is $F$-finite. More precisely, if $R=k[x_1,\ldots,x_n]$ is
a polynomial ring over a field $k$ of positive characteristic, then
every ideal $I$ in $R$ can be expressed as a generalized test ideal.

To see this, write $I=(z_1,\ldots,z_{\mu})$, and let $k_0$ be the
subfield of $k$ generated over the prime field $\bbF_p$ by the
coefficients of $z_1,\ldots,z_{\mu}$. Since $k_0$ is an extension of
finite type of a perfect field, it follows that $k_0$ is $F$-finite.
Therefore $S=k_0[x_1,\ldots,x_n]$ is also $F$-finite, and we may
apply Theorem~\ref{thm2} for $S$ to find $f\in S$ and $c\in\bbQ$
such that $\tau((fS)^c)=(z_1,\ldots,z_{\mu})S$. Since $R$ is free
over $S$, one can easily see that $\tau((fS)^c)R=\tau((fR)^c)$,
hence $I=\tau((fR)^c)$.
\end{remark}

It would be interesting to determine also in the singular case those
ideals that can be written as generalized test ideals. We end this
section with the following question of Shunsuke Takagi.

\begin{quest}\label{q-takagi}
Is the analogue of Theorem~\ref{thm2} true if we only assume that
the ring is strongly $F$-regular ?
\end{quest}

\section{Miscellaneous examples}

In this section we give several examples to show that the analogues
of several basic properties of multiplier ideals (which follow
easily from definition) fail for test ideals. We start by describing
the questions we will consider.

\begin{quest} \label{questions}
Let $(R,\frm)$ be an $F$-finite regular local ring of characteristic
$p>0$ with $d=\dim R \ge 1$. Let $f$ be a nonzero element of $R$,
and set $c=\fpt(f)$. Given $t > 0$, we put
$\tau(f^{t-})=\tau(f^{t-\varepsilon})$ for $0 < \varepsilon \ll 1$
(note that this is well-defined, since the $F$-jumping exponents of
$f$ are discrete; see \cite{BMS1}).
\begin{enumerate}
 \item[1)] Is the ideal $\tau(f^c)$ radical ?
 \item[2)] Suppose that $\tau(f^c)$ is $\frm$-primary.
Is there an $\frm$-primary ideal $\frb$ such that $f\in \frb$ and
$\fpt(f)=\fpt(\frb)$ ?
\item[3)] Does the inclusion
\[
 \frb^m \cdot \tau(f^{t-}) \cap \tau(f^t) \subseteq \frb^{m-d} \cdot \tau(f^t)
\]
hold for every $m\geq d$ and every $t>0$ ?
\item[4)] Does the analogue of the Generic Restriction Theorem for
multiplier ideals (see Theorem~\ref{thm30} below) hold for
generalized test ideals ?
\end{enumerate}
\end{quest}

We recall the argument for 1) and 2) in the case of multiplier
ideals. Suppose that $\fra$ is a nonzero ideal sheaf on the
nonsingular variety $X$ (over an algebraically closed field of
characteristic zero). Let $\pi\colon \widetilde{X}\to X$ be a log
resolution of the pair $(X,V(\fra))$. If
$\fra\OO_{\widetilde{X}}=\OO(-F)$, we write
\[
F=\sum_{i=1}^ra_iE_i,\qquad K_{\widetilde{X}/X}=\sum_{i=1}^rk_iE_i.
\]
Suppose that $c=\lct(\fra)$, hence $c=\min_i\frac{k_i+1}{a_i}$.

The analogue of 1) above holds since $\JJ(\fra^c)$ is the radical
ideal corresponding to $\cup_if(E_i)$, the union being over those
$i$ such that $c=\frac{k_i+1}{a_i}$. Moreover, suppose that $x\in X$
is a closed point corresponding to the ideal $\frm$. If
$\JJ(\fra^c)$ is $\frm$-primary, it follows that there is a divisor
$E_i$ lying over $x$, such that $c=\frac{k_i+1}{a_i}$. In this case,
for every $\ell>a_i$, we have ${\rm ord}_{E_i}(f)={\rm
ord}_{E_i}((f)+\frm^{\ell})$. Therefore $c\geq
\lct((f)+\frm^{\ell})$, and we get the assertion in 2), since the
reverse inequality is trivial.

The motivation for the question in 3) comes from its relevance to
uniform Artin-Rees results. The corresponding statement for
multiplier ideals is Theorem~3.1 in \cite{ELSV}. The proof uses only
the definition via log resolutions and Skoda's Theorem (which also
holds in the setting of test ideals). It is used to give an
effective uniform Artin-Rees statement for every ideal that can be
written as a multiplier ideal. Therefore, in light of our
Theorem~\ref{thm2}, a positive answer to 3) would have had very
strong consequences. It is conceivable that some weaker version of
3) might still hold, enough to give effective uniform Artin-Rees for
\textit{every} ideal in positive characteristic.

\bigskip

Our main source of counterexamples to the above questions is the
following proposition, giving a formula for all the test ideals of a
certain class of principal ideals.

\begin{prop} \label{concrete}
Let $p$ be a prime number, $n$ a positive integer, and let $R =
\bbF_p[[x_0,x_1,\ldots,x_n]]$ be a formal power series ring over
$\bbF_p =\bbZ/p\bbZ$. For any nonnegative integers
$\ell_1,\ldots,\ell_n$, we set
\[
 f = x_0^p+x_1^{\ell_1 p+1} + \cdots + x_n^{\ell_n p+1}
\quad \text{and} \quad I = (x_0,x_1^{\ell_1},\ldots,x_n^{\ell_n}).
\]
Then
\[
\tau(f^t) = \left\{
\begin{array}{cl}
R, & (0 \le t < \frac{~1~}{p}); \\[2mm]
I, & (\frac{~1~}{p} \le t < \frac{~2~}{p}); \\[2mm]
\vdots & \vdots \\[2mm]
I^{p-1}, & (\frac{p-1}{p} \le t < 1); \\[2mm]
fR, & (t=1).
\end{array}\right.
\]
In particular,
\begin{enumerate}
 \item $\fpt(f) = \frac{~1~}{p}$ and $\tau(f^{\fpt(f)}) = I$.
 \item For every $t \in \bbR_{\ge 0}$, we have
\[
 \tau(f^{t}) = f^{\lfloor t \rfloor} \,
I^{\lfloor p (t-\lfloor t \rfloor) \rfloor},
\]
where $\lfloor \alpha \rfloor$ denotes the largest integer
$\leq\alpha$.
 \item The set of $F$-jumping exponents of $f$ is $\frac{~1~}{p}\,\bbZ_{\ge 0}$.
\end{enumerate}
\end{prop}

\begin{proof}

It is enough to show that $\tau(f^t) = I^r$ for
$t\in\left[\frac{r}{p},\frac{r+1}{p}\right)$ and for every
$r=0,1,\ldots,p-1$. The other assertions follow from this and
Skoda's Theorem. First, we show the following
\begin{flushleft}
 {\bf Claim 1:} $\tau(f^{r/p}) = I^r$.
\end{flushleft}
Since we have
\begin{eqnarray*}
f^{\lceil (r/p)p^e \rceil} &= & f^{rp^{e-1}} =
\left(x_0^{p^e}+x_1^{\ell_1p^e+p^{e-1}}+\cdots
+ x_n^{\ell_np^e+p^{e-1}}\right)^r \\
& = & \sum_{\stackrel{\scriptstyle i_0,\ldots,i_n \ge 0}{i_0 +
\cdots + i_n = r}} \frac{r!}{i_0!i_1!\cdots i_n!}\; \left(x_0^{i_0}
x_1^{\ell_1i_1} \cdots x_n^{\ell_n i_n} \right)^{p^e} \;
x_1^{i_1p^{e-1}}\cdots x_n^{i_np^{e-1}}
\end{eqnarray*}
and since $\left\{\frac{r!}{i_0!i_1!\cdots i_n!}\;
x_1^{i_1p^{e-1}}\cdots x_n^{i_np^{e-1}}\right\}$ is part of a free
basis of $R$ over $R^{p^e}$, we obtain that
\[
\left(f^{\lceil(r/p)p^e \rceil} \right)^{[1/p^e]}=
(x_0,x_1^{\ell_1},\ldots,x_n^{\ell_n})^r.
\]
Since this holds for every $e\geq 1$, we get our claim.

\par \vspace{2mm}
In order to prove that $\tau(f^t) = I^r$ when $\frac{r}{~p~} < t <
\frac{r+1}{~p~}$, we put $t=\frac{r+1}{~p~}-\varepsilon$, $0 <
\varepsilon < \frac{1}{~p~}$. It follows from Claim~1 that it is
enough to show that $I^r\subseteq\tau(f^t)$. We fix a sufficiently
large integer $e$ such that  $s:=\lfloor \varepsilon p^e \rfloor \ge
1$. We have
\begin{eqnarray*}
 f^{\lceil tp^e \rceil}
&=& \left(x_0^p+x_1^{\ell_1p+1} +
\cdots + x_n^{\ell_n p+1} \right)^{(r+1)p^{e-1}-s} \\
&=& \sum_{\stackrel{\scriptstyle a_0,\ldots,a_n \ge 0}{a_0 + \cdots
+ a_n =(r+1)p^{e-1} -s}}
\!\!\!\!\!\frac{((r+1)p^{e-1}-s)!}{a_0!a_1!\cdots a_n!} \;\;
x_0^{pa_0} x_1^{(\ell_1p+1)a_1}\cdots x_n^{(\ell_np+1)a_n}.
\end{eqnarray*}
\par
In order to complete the proof, it is enough to show that for every
$(n+1)$-tuple of nonnegative integers $(i_0,i_1,\ldots,i_n)$ such
that $i_0+i_1+\cdots + i_n =r$, we have
\[
y:=x_0^{i_0}x_1^{\ell_1i_1}\cdots x_n^{\ell_ni_n}\in \left(f^{\lceil
tp^e\rceil}\right)^{[1/p^e]}.
\]
If we put $a_0 = (i_0+1)p^{e-1}-s$, $a_j = i_jp^{e-1}$ for
$j=1,\ldots,n$, then we have
\[
 a_0,a_1,\ldots, a_n \ge 0, \quad a_0+a_1+\cdots + a_n =
 (r+1)p^{e-1}-s
\]
and
\[
 x_0^{pa_0}x_1^{(\ell_1p+1)a_1}\cdots x_n^{(\ell_np+1)a_n}
= \left(x_0^{i_0}x_1^{\ell_1i_1}\cdots x_n^{\ell_ni_n} \right)^{p^e}
\; x_0^{p^e-sp} x_1^{i_1p^{e-1}}\cdots x_n^{i_np^{e-1}}.
\]
Therefore it is enough to prove the claim below. Note that the claim
implies that $f^{\lceil tp^e\rceil}$ can be written as $y_1^{p^e}
g_1 + \cdots + y_{\mu}^{p^e} g_{\mu}$, such that
$I^r=(y_1,\ldots,y_{\mu})$ and $\{g_1,\ldots,g_{\mu}\}$ is part of a
free basis of $R$ over $R^{p^e}$.

\begin{flushleft}
{\bf Claim 2:}
\begin{enumerate}
 \item $\frac{((r+1)p^{e-1}-s)!}{a_0!a_1!\cdots a_n!}  \not \equiv 0 \pmod{p}$.
\vspace{2mm}
 \item Let $b_0,b_1,\ldots,b_n \ge 0$ be integers
with $b_0+b_1 + \cdots + b_n = (r+1)p^{e-1} -s$. If there exist
$t_0,t_1,\ldots,t_n \in \bbZ$ such that
\[
 pb_0 - pa_0 = t_0p^e,\qquad (\ell_j p+1)(b_j-a_j) = t_j p^e \; (j=1,\ldots,n),
\]
then $b_0=a_0$, $b_1=a_1,\ldots, b_n=a_n$.
\end{enumerate}
\end{flushleft}
\par
In order to prove (1), we use the fact that for every integer $N$,
the order of $p$ in $N!$ is $\sum_{m\geq 1}\lfloor N/p^m\rfloor$.
Note that if $1\leq m\leq e-1$, then we have
\[
\lfloor(a_0+a_1+\cdots a_n)/p^m\rfloor =\lfloor
a_0/p^m\rfloor+\sum_{j=1}^ni_jp^{e-1-m}=\sum_{j=0}^n\lfloor
a_j/p^m\rfloor.
\]
On the other hand, $a_0+a_1+\cdots+a_n<p^e$. This shows that the
order of $p$ in $\frac{((r+1)p^{e-1}-s)!}{a_0!a_1!\cdots a_n!}$ is
zero.

\par \vspace{2mm}
We now prove (2).  Since ${\rm gcd}(p,\ell_j p+1)=1$, we have
$p^e\mid (b_j-a_j)$ for every $1\leq j\leq n$. Therefore we can
write $b_j-a_j = u_jp^e$ for every $j$ as above, and suitable
$u_j\in\bbZ$. Using $b_j = (i_j+pu_j)p^{e-1} \ge 0$, we deduce $i_j
+ pu_j \ge 0$, hence $u_j \ge 0$ (recall that
 $i_0+\cdots + i_n =r < p$). On the
other hand, since $b_0 = (i_0+1+t_0)p^{e-1}-s \ge 0$, we get
$i_0+1+t_0 > 0$ and thus $t_0 \ge -i_0 > -p$. Moreover, $a_0+\cdots
+ a_n = b_0 + \cdots + b_n$ yields $(u_1+\cdots + u_n)p+t_0 =0$.
Therefore $a_j=b_j$ for every all $j$. This completes the proof of
Claim~2, and also the proof of the proposition.
\end{proof}

\begin{exam} \label{counterex}
Let $R = \bbF_2\llbracket x,y,z\rrbracket$, $f = x^2+y^5+z^5$, and
put $\fra_N = (f) +(x,y,z)^N$ for every $N \ge 1$.
\begin{enumerate}
 \item $\fpt(f) = \frac{1}{~2~}$.
 \item $\tau(f^{\fpt(f)}) = (x,y^2,z^2)$ is an $\frm$-primary ideal,
but it is \textit{not} radical (hence this gives a counterexample to
1) in Question~\ref{questions}.
 \item $\fpt(\fra_N) > \fpt(f)=\frac{1}{~2~}$ for every $N \ge 1$
 (hence this gives a counterexample to 2) in Question~\ref{questions}).
\end{enumerate}
\end{exam}

\begin{proof}
(1) and (2) follow from Proposition~\ref{concrete}. In order to see
that (3) indeed says that we get a counterexample to 2) in
Question~\ref{questions}, note that if $\frb$ is an $\frm$-primary
ideal containing $f$, then there is $N\geq 1$ such that $\fra_N
\subseteq\frb$. Hence $\fpt(\frb)\geq\fpt(\fra_N)>\fpt(f)$.
\par
It is enough to prove the assertion in (3) for every $N=2^{e-2}$,
where $e\geq 5$. We show that in this case $\fra_N^{2^{e-1}}\not
\subseteq (x^{2^e},y^{2^e}, z^{2^e})$, hence $\tau(\fra_N^{1/2})=R$.
Consider
\[
h:= f^{2^{e-1}-4}x^Ny^Nz^{2N}\in \fra_N^{2^{e-1}}.
\]
If $a=2(2^{e-1}-4-2^{e-3})+2^{e-2}=2^e-8$, $b=5\cdot
2^{e-3}+2^{e-2}=7\cdot 2^{e-3}$, and $c=2^{e-1}$, then the monomial
$x^ay^bz^c$ is not in $(x^{2^e},y^{2^e}, z^{2^e})$, and its
coefficient in $h$ is ${{2^{e-1}-4}\choose 2^{e-3}}$. In order to
show that this coefficient is nonzero, we compute the order of $2$
in ${{2^{e-1}-4}\choose 2^{e-3}}$. This order is equal to
\[
\sum_{i=1}^{e-2}\left(\lfloor (2^{e-1}-4)/2^i\rfloor -\lfloor
2^{e-3}/2^i\rfloor-\lfloor (2^{e-1}-4-2^{e-3})/2^i\rfloor\right)
\]
\[
=\lfloor(2^{e-1}-4)/2^{e-2}\rfloor-\lfloor
(2^{e-1}-4-2^{e-3})/2^{e-2}\rfloor=1-1=0.
\]
This concludes the proof of (3).
\end{proof}

\begin{remark} \label{Schwede}
Karl Schwede \cite{Sch} has recently introduced the notion of
\textit{sharp $F$-purity}. He proved that
 if $c= \fpt(f) < 1$ is such that the denominator of $c$ is not
 divisible by $p$, then the ideal $\tau(f^c)$ is radical;
see Corollary 4.3 and Remark 5.5 in \emph{loc. cit}. It would be
very interesting to see whether assuming that the denominators of
the jumping numbers of $f$ are not divisible by $p$ would imply
other good properties of the generalized test ideals of $f$.
\end{remark}

We consider now the third problem in Question~\ref{questions}.

\begin{exam}\label{AR}
Let $p$ be a prime, $R=\bbF_p\llbracket x,y\rrbracket$ and
$f=x^p+y^{\ell p+1}$, for some $\ell\geq 3$. It follows from
Proposition~\ref{concrete} that $\fpt(f)=1/p$ and
$\tau(f^{1/p})=(x,y^{\ell})$. If we take $\frb=(y)$ and $t=1/p$,
then we see that
\[
\frb^{\ell}\cdot\tau(f^{t-}) \cap \tau(f^{t}) =\frb^{\ell}\cap
(x,y^{\ell})=(y^{\ell})\not\subseteq
\frb^{\ell-2}\cdot\tau(f^t)=(y^{\ell-2})\cdot (x,y^{\ell}),
\]
giving thus a counterexample to 3) in Question~\ref{questions}.
\end{exam}

\bigskip

We conclude this section with a discussion of the analogue of the
Generic Restriction Theorem for multiplier ideals in the
characteristic $p$ setting. Let us recall the result in
characteristic zero (see \cite{Laz}, Theorem~9.5.35 and
Example~9.5.37).

\begin{thm}\label{thm30}
Let $f\colon X\to S$ be a smooth surjective morphism of nonsingular
complex algebraic varieties. If $\fra$ is a sheaf of ideals on $X$,
then there is an open subset $U\subseteq S$ such that
$$\mathcal{J}(X,\fra^c)\cdot\mathcal{O}_{X_s}
=\mathcal{J}(X_s,(\fra\cdot\mathcal{O}_{X_s})^c)$$ for every $s\in
U$ and every positive $c$ {\rm (}here $X_s$ denotes the fiber
$f^{-1}(s)${\rm )}.
\end{thm}

We show now that the analogue of this result fails for test ideals.
Suppose, for simplicity, that $k$ is an algebraically closed field
of positive characteristic, and consider $f\in
R=k[x_1,\ldots,x_n,y]$. Let us denote by $\{u_j\}_j$ the monomials
$x_1^{a_1}\cdots x_n^{a_n}$, where $0\leq a_i\leq p-1$ for every
$i$. We write
\begin{equation}\label{eq1}
f=\sum_{i=0}^{p-1}y^i\sum_ju_jg_{ij}(x,y)^p,
\end{equation}
for some $g_{ij}\in R$. Arguing as in the proof of
Proposition~\ref{Precise}, we see that
\begin{equation}
\tau(f^{1/p})=(f)^{[1/p]}= (g_{ij}(x,y)\mid i,j).
\end{equation}

\par
On the other hand, let us put $f_{\lambda}(x):=f(x,\lambda)\in
k[x_1,\ldots,x_n]$ for every $\lambda\in k$. Note that we have
\begin{equation}\label{eq2}
f_{\lambda}=\sum_j u_j\sum_{i=0}^{p-1}g_{ij}(x,\lambda)^p\lambda^i,
\end{equation}
hence we deduce
\begin{equation}
\tau(f_{\lambda}^{1/p})=(f_{\lambda})^{[1/p]}=
\left(\sum_{i=0}^{p-1}\lambda^{i/p}g_{ij}(x,\lambda)\mid j\right).
\end{equation}

\begin{exam}\label{semicont1}
Consider $f\in k[x_1,x_2,y]$ given by $f(x_1,x_2,y)=x_1^p+x_2^py$.
The above discussion implies that $\tau(f^{1/p})=(x_1,x_2)$, while
for every $\lambda\in k$ we have
$\tau(f_{\lambda}^{1/p})=(x_1+\lambda^{1/p}x_2)$. This gives a
negative answer to 4) in Question~\ref{questions}.
\end{exam}

The main application of Theorem~\ref{thm30} is to prove the
semicontinuity of log canonical thresholds. In spite of the above
example, we will see in the next section that the analogous result
for $F$-pure thresholds holds.

\section{Semicontinuity of $F$-pure thresholds. }

The following theorem is the analogue of the Semicontinuity Theorem
for log canonical thresholds (see \cite{Laz}, Example~9.5.41).

\begin{thm}\label{thm3}
Let $f\colon R\to S$ be an algebra homomorphism between two
$k$-algebras of finite type, where $k$ is a field of characteristic
$p$, with $[k\colon k^p]<\infty$. We assume that all fibers of $f$
are nonsingular, of pure dimension $d$. Let $\phi\colon S\to R$ be a
ring homomorphism such that $\phi\circ f= {\rm id}_R$, and for every
$\frq\in {\rm Spec}(R)$, we put $\frq'=\phi^{-1}(\frq)$. For every
ideal $\fra$ in $S$ such that $\fra\subseteq\frq'$ for all $\frq\in
{\rm Spec}(R)$, and for every nonnegative $c$, the set
$$\{\frq\in {\rm Spec}(R)\mid \fpt(\fra \cdot S_{\frq'}/q S_{\frq'})\geq c\}$$
is open in ${\rm Spec}(R)$.
\end{thm}

\begin{proof}
Note that for every $\frq\in {\rm Spec}(R)$ we have $[k(\frq)\colon
k(\frq)^p]<\infty$, hence the ring $S_{\frq'}/\frq S_{\frq'}$ is
$F$-finite and regular.
 Consider a surjective morphism of $R$-algebras $g\colon
T=R[x_1,\ldots,x_n]\to S$. We claim that we may replace $S$ by
$R[x_1,\ldots,x_n]$. Indeed, it follows from Proposition~3.6 in
\cite{BMS1} that if we write $\fra=\frb/{\rm ker}(g)$ and
$\frq'=\frq''/{\rm ker}(g)$, then
$$\fpt (\fra \cdot S_{\frq'}/\frq S_{\frq'})+n-d=\fpt(\frb \cdot T_{\frq''}/\frq
T_{\frq''}).$$ This proves our claim. Moreover, note that if
$\phi\colon S=R[x_1,\ldots,x_n]\to R$ is given by $\phi(x_i)=b_i$,
then we may consider the automorphism of $R$-algebras $\rho\colon
S\to S$ given by $\rho(x_i)=x_i+b_i$. After replacing $\fra$ by
$\rho(\fra)$,  we may also assume that $\phi(x_i)=0$ for every $i$.
We see that for every $\frq\in {\rm Spec}(R)$, we are interested in
the $F$-pure threshold of $\fra\cdot
k(\frq)[x_1,\ldots,x_n]_{(x_1,\ldots,x_n)}$, that we denote by
$\fpt_0\left(\fra\cdot k(\frq)[x_1,\ldots,x_n]\right)$.

Let us choose generators $g_1,\ldots,g_m$ for $\fra$, and let
$D=\max_i\{\deg(g_i)\}$. It follows from Proposition~3.8 in
\cite{BMS1} that there is $N=N(D,n,m)$ such that the denominator of
every $F$-jumping exponent of an ideal of the form $\fra \cdot
k(\frq)[x_1,\ldots,x_n]$ (for $\frq\in {\rm Spec}(R)$) is $\leq N$.
Note that $\fpt_0\left(\fra\cdot k(\frq)[x_1,\ldots,x_n]\right)$ is
an $F$-jumping exponent of $\fra\cdot k(\frq)[x_1,\ldots,x_n]$
(though it might be larger than the $F$-pure threshold of this
ideal). Using also the fact that the $F$-pure threshold of an ideal
in a regular ring of dimension $n$ is $\leq n$, we deduce that the
set
\[
\{\fpt_0(\fra \cdot k(\frq)[x_1,\ldots,x_n])\mid \frq\in {\rm
Spec}(R)\}
\]
is finite.

In particular, in order to prove the theorem, we may choose the
largest element $c'$ in the above set, with $c'<c$. It is enough to
show that the set
$$A_{c'}:=\{\frq\in {\rm Spec}(R)\mid \fpt_0(\fra\cdot
k(\frq)[x_1,\ldots,x_n])\leq c'\}$$ is closed. Using the description
of the $F$-pure threshold in (\ref{fpt}) in \S 2, we see that
$A_{c'}=\cap_{e\geq 1}A_{c',e}$, where
$$A_{c',e}=\{\frq\mid \fra^{\lfloor c'p^e\rfloor+1}\subseteq
(x_1^{p^e},\ldots,x_n^{p^e})\,{\rm in}\,k(\frq)[x_1,\ldots,x_n]\}.$$
Note that if we consider all $g^{\ell}:=g_1^{\ell_1}\cdots
g_m^{\ell_m}$, with $\sum_{i}\ell_i=\lfloor c'p^e\rfloor+1$, then
$A_{c',e}$ is the set of primes $\frq$ containing all the
coefficients of monomials not in $(x_1^{p^e},\ldots,x_n^{p^e})$, in
all $g^{\ell}$ as above. Therefore each $A_{c',e}$ is a closed
subset of ${\rm Spec}(R)$, hence $A_{c'}$ is closed, too.
\end{proof}

\medskip

{\bf Acknowledgements}. We are indebted to Shunsuke Takagi and to
Kei-ichi Watanabe for inspiring discussions.

\end{document}